\documentclass[12pt,fleqn,legno]{article}
\usepackage{amsmath}
\usepackage{amsfonts}
\usepackage{amsthm}

\newcommand{\Rset}{\mathbb{R}}
\newcommand{\Cset}{\mathbb{C}}

\newtheorem{thm}{THEOREM}

\newtheorem{prop}{PROPOSITION}
\newtheorem{cor}{COROLLARY}

\theoremstyle{remark}
\newtheorem{rem}{REMARK}

\begin{document}

\title{On Relations Between Urbanik and Mehler
Semigroups\footnote{Research funded by a grant MEN Nr 1P03A04629,
2005-2008.}}

\author{Zbigniew J. Jurek (Wroclaw)}

\date{November 14, 2008}

\maketitle

\begin{quote} \textbf{Abstract.} It is shown that
operator-selfdecomposable measures, or more precisely their Urbanik
decomposability semigroups, induce generalized Mehler semigroups of
bounded linear operators. Moreover, those semigroups can be
represented as random integrals of operator valued functions with
respect to stochastic L\'evy processes. Our Banach space setting is
in the contrast with the Hilbert spaces on which so far and most
often the generalized Mehler semigroups were studied. Furthermore,
we give new proofs of the random integral representation.


\end{quote}

\medskip
\emph{An abbreviated title:} Urbanik and Mehler semigroups.

\newpage
The theory of operator limiting distributions in probability theory
has its origin in 70's of the last century. Its development on the
Euclidean spaces was summarized in the monograph by Jurek and Mason
(1993). That theory is based on the principle that normalization of
random variables with values in $E$ (or sums of those random
variables) should be consistent with the structure of the state
space $E$. Thus the linear operators are the proper normalization
for the linear spaces or normalization by the group automorphisms
for the case of group valued variables .

The most important new tool in that setting is \emph{the operator
decomposability semigroup} introduced by K.Urbanik in 1972. Namely,
with a probability measure $\mu$ one associates a family
$\textbf{D}(\mu)$ of linear bounded operators $A$ that "divide"
$\mu$ in  a sense that $\mu=A\mu\ast\nu_A$, for some another
probability measures $\nu_A$, i.e.,
\[
\textbf{D}(\mu)=\{A: \mu=A\mu\ast\nu_A \ \mbox{for some} \ \nu_A\}
\]
Note that operators $0$, $I$ are always in $\textbf{D}(\mu)$  and
that it is indeed a semigroups. If the Urbanik semigroup
$\textbf{D}(\mu)$ contains an one-parameter semigroup $T_t, t\ge 0$,
and one defines $\rho_t:=\nu_{T_t}$ then, by an iteration, one
arrives to the equation
\begin{equation}
\ \ \ \  \rho_{t+s}=\rho_t\ast T_t\rho_s,  \ \ \mbox{for all} \ \ s,
t \ge 0, \ \ \ \
\end{equation}
provided a cancelation is permitted; cf. Jurek (1982), Jurek-Vervaat
(1983). Such convolution equations also were called
\emph{measure-valued cocycles}; cf. Jurek-Hofmann(1996). [Note that
the Grothendick type diagram, on p. 755 there, is still not
completed !]

On the other hand, if operators $\mathcal{T}_t$, given by
\[
(\mathcal{T}_tf)(x)=\int_Ef(T_tx+z)\rho_t(dz), \ f\in C_b(E) ,\
t\ge 0,
\]
define an one-parameter semigroup on $C_b(E)$ then the measures
$\rho_t$'s must satisfy the cocycle relation (1). The above families
of operators are called \emph{generalized Mehler semigroups}, for
short: Mehler semigroups. Cf. for instance Bogachev, Roeckner and
Schmuland (1996) and references therein. However, one should be
aware that they worked on Hilbert spaces and had more restrictive
assumption (differentiability) on the Fourier transforms of $\rho_t,
t\ge0$. On the other hand, cocycle equations (1) were considered on
non-linear structures like spaces of measures; cf. Li (2002).

The main result here is that on a Banach space under the continuity
of the mapping  $t\to \rho_t$ in many cases Mehler semigroups are of
the form
\[
(\mathcal{T}_t\,f)(x) = \mathbf{E}
[f(T_tx+\int_{(0,t]}T_{t-s}\,dY(s))], \ f\in C_b(E),
\]
for some stochastic L\'evy process $Y$; cf. below for the details.

\medskip
In all papers dealing with the operator limit distributions (for
instance: Jurek (1982), (1983), (1988), Jurek-Vervaat (1983) and
others) the primary goal was the random integral representation
(RIR) of  measures $\mu$, whose Urbanik semigroups $\textbf{D}(\mu)$
contain one-parameter semigroup of operators continuous in the
operator norm topology. Consequently, the solutions to the cocycles
equations (1) were auxiliary steps in the main proofs and might have
been overlooked.

These two subjects, i.e., Urbanik and Mehler semigroups, seem to be
developed independently of each other, although Chojnowska-Michalik
(1987) mentioned operator limit distributions theory. See also the
acknowledgment at the end of this paper.

Our aim here is to show how the operator limit distributions theory
and its techniques (like the random integral method) can produce new
results and proofs in the theory of Mehler semigroups on Banach
spaces; cf. in particular Proposition 4, Theorem 1 and Corollaries 2
and 3 below.

Last but not least, let us stress again that  our presentation here
is in the generality of Banach spaces while Hilbert space setting
was often the case for the generalized Mehler semigroups.

\medskip
\textbf{1. Basic notions and notations.} Let $E$ denotes \emph{a
real separable} Banach space with a norm $||.||$, let $\bf{E}$nd(E),
or simply $\bf{E}$nd,  denotes the algebra of \emph{all} bounded
linear operators on $E$. In $\bf{E}$nd(E) we take \emph{the strong
operator topology},  i.e., $A_n \stackrel{s}{\to}A$ means that
\emph{for each} $x\in E$, \ \ \ $\lim_{n\to \infty}||A_n\,x-Ax||=
0$. Of a particular interest are the \emph{$C_0$-one-parameter
semigroups} $\textbf{T}=(T_t, t\ge 0)$ in \textbf{E}nd; that is we
have: $ T_0=I, \ \ T_t(T_s\,x)=T_{t+s}x, \ \ \mbox{for } \ t,s\ge0,
\ x\in E \ \ \mbox{and for each $x$ the functions}\ \ t\to T_tx \ \
\mbox{are continuous.}$

Let $\mathcal{P}(E)$ or just $\mathcal{P}$, denotes the family of
\emph{all Borel probability measures} on $E$ endowed with the
\emph{convolution} $"\ast"$ operation and the \emph{weak
convergence} topology, in symbols: $"\Rightarrow"$. Thus
\[
\mu_n \Rightarrow \mu \ \ \mbox{iff} \ \
\int_Ef(x)\mu_n(dx)\to\int_Ef(x)\mu(dx), \ \ \mbox{for each} \ \
f\in C_b(E);
\]
\noindent where  $C_b(E)$ stands for real-valued continuous bounded
functions on $E$; (weak*-topology in $C_b(E)$). For the probability
theory on Banach spaces see Araujo-Gine (1980) or Linde (1986).

Let $E'$ be the topological dual Banach space and let $<.,.>$
denotes the bilinear form between $E'$ and $E$. Recall that for a
measure $\mu$ or  an E-valued random variable $\xi$ with probability
distribution $\mu$, the function
\begin{equation*}
\widehat{\mu} : E'\to \Cset \ \  \mbox{giveny by} \ \ \ \
\widehat{\mu}(y):=\int_E\,e^{i<y,x>}\,\mu(dx)=\textbf{E}[e^{i<y,\xi>}],
\end{equation*}
is called the \emph{Fourier transform} (or the \emph{characteristic
function}) and that it uniquely determines the measure $\mu$; above
$\textbf{E}[.]$ denotes the expectation operator.

Finally, for $A\in \textbf{E}nd$ and $\mu \in \mathcal{P}$ we define
$A\mu \in \mathcal{P}$, the image of $\mu$ through a mapping $A$, as
follows:
\[
\ \ \ \ \ (A\mu)(\mathcal{E}):\,=\,\mu(\{x\in E:
Ax\in\mathcal{E}\})\ \ \mbox{for all Borel subsets $\mathcal{E}$ of
E}.
\]
Equivalently, in terms of integrals, it means that
\[
\int_Ef(x)(A\mu)(dx)=\int_Ef(Ax)\mu(dx), \ \ \mbox{for all} \ \ f\in
C_b(E).
\]
\noindent In other words, if $\xi$ is an E-valued random variable
with probability distribution $\mu$ then the random variable $A\xi$
has probability distribution $A\mu$. Having that in mind we
immediately get the equalities
\begin{equation}
A(\mu\ast\nu)=A\mu \ast A\nu, \ \ \ A(B\mu)= (AB)\mu , \
\widehat{(A\mu)}(y)= \hat{\mu}(A^*y), y\in E^{\prime},
\end{equation}
for all linear bounded operators $A, B\in \textbf{E}nd$ and all
measures $\mu, \nu \in \mathcal{P}$.

Finally for the future reference let us quote here that
\begin{equation}
\mbox{if} \ A_n\stackrel{s}{\to} A \ \mbox{and} \ \mu_n \Rightarrow
\mu \ \mbox{ then} \ \  A_n \mu_n \Rightarrow A\mu.
\end{equation}
Proof can be found in Jurek (1983), Proposition 1.1 or in
Jurek-Mason (1993), Proposition 1.7.2, on p. 24.

\medskip
 \textbf{2. The Urbanik decomposability semigroups.}
With $\mu\in\mathcal{P}$ we associate its \emph{Urbanik
decomposability semigroup} $\textbf{D}(\mu)$ defined as follows
\begin{equation}
\textbf{D}(\mu):=\{A\in \textbf{E}nd: \mu= A\mu \ast \nu_A, \
\mbox{for some} \ \nu_A\in \mathcal{P}\}
\end{equation}
Obviously, the linear operators $0$ (zero) and $I$ (identity) are in
all $\textbf{D}(\mu)$ with $\nu_0=\mu$ and $\nu_I=\delta_0$ in (4)
and the semigroup property, under composition of operators, follows
from (2). It is interesting that some purely probabilistic
properties  of $\mu$ are equivalent with some algebraic and
topological properties of its Urbanik $\textbf{D}(\mu)$
decomposability semigroup; cf. Urbanik (1972), (1978), Jurek-Mason
(1993)

In the operator-limit distribution theory the operator topology is
used in  $\textbf{D}(\mu)$. However, even for the strong operator
topology we also have that
\begin{prop}
(i) \ The Urbanik decomposability semigroups $\textbf{D}(\mu)$, in
$\mathbf{E}nd$, are closed in the strong operator topology.

(ii) \ If $\hat{\mu}(y)\neq 0$ for all $y\in E'$, $A_n\in
\textbf{D}(\mu)$ and $A_n\stackrel{s}{\to}A$ then $A\in
\textbf{D}(\mu)$ and, in (4), we have that \ $\nu_{A_n}\Rightarrow
\nu_A$.

(iii) \ If $\mu=A_n\mu\ast\nu_{A_n}$ and $A_n \stackrel{s}{\to} 0$
then $\nu_{A_n}\Rightarrow \mu$.
\end{prop}
\emph{Proof.} (i) For $A_n\in \textbf{D}(\mu)$ we have
\begin{equation}
\mu=A_n\mu\ast\nu_{A_n} \ \ \mbox{for some } \ \ \nu_{A_n}\in
\mathcal{P}.
\end{equation}
Further, if $A_n\stackrel{s}{\to}A$ then, by (3), $A_n\mu\Rightarrow
A\mu$. Consequently, $\{\nu_{A_n}, n=1,2,...\} \subset\mathcal{P}$
is conditionally compact (uniformly tight); cf. Parthasarathy
(1967), Chapter III, Theorem 2.1 or Jurek-Mason (1993), Theorem
1.7.1, Thus passing to a subsequence in (5) we get\
$\mu=A\mu\ast\nu$ \ for some accumulation point $\nu$ of the
sequence $(\nu_{A_n}$, n=1,2,... ). Consequently, $A\in
\textbf{D}(\mu)$, which proves (i).

\noindent (ii) From (i) we get that $A\in \textbf{D}(\mu)$. Since
$\hat{\mu}(A_n^{*}y)\neq 0$ ($A^{*}$ is the conjugate bounded linear
operator), from (5), we infer that
$\lim_{n\to\infty}\widehat{\nu}_{A_n}(y)$ exists. This and the
conditional compactness of $(\nu_{A_n}, n=1,2,...)$ implies the weak
convergence $\nu_{A_n}\Rightarrow \nu_A$ in (ii).

\noindent (iii) Simply note, by  (3), that $A_n\mu\Rightarrow
\delta_0$ and thus $\hat{\mu}(A_n^{*}y)\to 1$ for all $y\in E'$.
Hence, as in the proof of (ii), we conclude $\nu_{A_n}\Rightarrow
\mu$. This completes the proof of Proposition 1.

\medskip
We will say that $\nu\in\mathcal{P}$ is  \emph{an operator
convolution factor} of $\mu$ if there exists $A\in\textbf{E}nd$ such
that $\mu=A\mu\ast\nu$. By $O\textbf{F}(\mu)$ we denote the totality
of the operator convolution factors of $\mu$.
\begin{prop}
Let $\mu\in\mathcal{P}$ be such that $\hat{\mu}(y)\neq 0$ for all $y\in
E'$. Then \ $O\textbf{F}(\mu)=\{\nu_A: A\in\textbf{D}(\mu)\}$ with
binary operation $\diamond$ given by $\nu_A \diamond
\nu_B:\,=\nu_A\ast\,A\nu_B$ is a non-commutative semigroup.
Moreover, $\nu_A \diamond \nu_B=\nu_{AB}$, $\nu_I\equiv\delta_0$ is
the neutral element and
$\nu_A\diamond\nu_0=\nu_0\diamond\nu_A=\nu_0\equiv\mu$.
\end{prop}
\emph{Proof.} For $A, B \in \textbf{D}(\mu)$ we have
\[
\mu = A\mu\ast \nu_A=A(B\mu\ast\nu_B)\ast\nu_A= (AB)\mu\ast
(A\nu_B\ast\nu_A)=(AB)\mu\ast\nu_{AB},
\]
because of (2) and the fact that $AB\in \textbf{D}(\mu)$ as well.
Hence  \ $\nu_A \diamond \nu_B=\nu_A\ast A\nu_B=\nu_{AB}$, \ because
of \ $\hat{\mu}(y)\neq0$. Furthermore,
\[
(\nu_A\diamond \nu_B)\diamond \nu_C=\nu_{AB}\diamond
\nu_C=\nu_{(AB)C}=\nu_A \diamond (\nu_B\diamond\nu_C),
\]
which proves the associativity of the operation $\diamond$. The rest
follows from the equalities $\nu_I\equiv\delta_0$ and
$\nu_0\equiv\mu$;\ comp. (4).

\begin{rem}
If $A_n\in \textbf{D}(\mu)$, $\hat{\mu}(y)\neq 0$ for all $y\in E'$
and $\nu_{A_n}\Rightarrow\rho_2\in\mathcal{P}$ then
$A_n\mu\Rightarrow\rho_1$, for some $\rho_1\in \mathcal{P}$ and
$\mu=\rho_1\ast\rho_2$, that is, $\rho_2$ is \emph{convolution
factor} of $\mu$. But is it \emph{an operator convolution factor} ?
Can we write $\rho_1=A\mu$, for some $A\in\textbf{E}nd$ ? In a case
of $E=\Rset^d$ (finite dimensional space) and \emph{full measures}
the answer is affirmative; see Lemma 2.2.9  and Section 2.5 in
Chapter II of Jurek-Mason (1993).
\end{rem}

We say that $\mu$ is \emph{operator-selfdecomposable} on $E$, in
symbols $\mu\in OS$, if the Urbanik semigroup $\textbf{D}(\mu)$
contains (at least one) $C_0$-semigroup $\textbf{T}=(T_t,t\ge0)$. \
When  a semigroup \textbf{T} is fixed then we write that $\mu\in
OS(\textbf{T})$ and say that $\mu$ \emph{is
\textbf{T}-decomposable}.

\begin{rem}
It is also important to realize that originally in Urbanik (1978)
there were operator continuous one-parameter semigroups
$T_t=\exp\,t\textbf{V}, t\ge 0$ such that $\lim_{t\to \infty}T_t=0$
(in the operator norm) . Furthermore, Urbanik primarily dealt with
limit distributions of sequences of partial sums of E-valued
variables normalized by arbitrary bounded linear operators. Similar
approach was taken in Jurek (1983), however with \emph{specified}
normalizing operators but with the strong operator topology.  For
the theory of operator-selfdecomposable (and operator-stable)
measures cf. Jurek-Mason (1993) and references therein.
\end{rem}
Explicitly, we have
\begin{equation}
\mu\in OS((T_t,t\ge0))\ \ \mbox{iff} \ \ \forall\, (t \ge 0)\,
\exists \, (\nu_{T_t}\in\mathcal{P}) \ \ \mu=T_t\mu\ast\nu_{T_t}
\end{equation}
\noindent Or equivalently, by Proposition 2, in terms of the
operator convolution factors semigroup $O\textbf{F}(\mu)$, we have
that
\begin{multline}
\mu\in OS(\,(T_t,t\ge0)\,)\ \ \mbox{iff} \ \
 \ \ \ \ \ \ \\ \exists\,(\{\rho_t, t\ge
0\}\subset(O\textbf{F}(\mu), \diamond))\ \ \forall(s,t\ge 0)\,\,
\rho_t\diamond\rho_s=\rho_{t+s},
\end{multline}
where $\rho_t:=\nu_{T_t}, \ t\ge 0$.

\medskip
\textbf{3. The Generalized Mehler semigroups.} For an operator $A\in
\textbf{E}nd(E)$ and a probability $\mu\in\mathcal{P}(E)$, let us
define the linear operator $\mathcal{A}^{(\mu)}$ as follows
\begin{equation}
\mathcal{A}^{(\mu)}:   C_b(E)\to C_b(E) \ \
(\mathcal{A}^{(\mu)}\,f)(x):\,=\int_E\,f(Ax+z)\mu(dz), \ x\in E.
\end{equation}
Note that $\mathcal{A}^{(\mu)}$ can viewed as the convolution of a
function $f$ with a measure $\delta_{Ax}\ast \mu$. Here are some
elementary properties of those operators.

\begin{prop}
(i) The operator \ $\mathcal{A}^{(\mu)}$ \ uniquely determines a
measure $\mu\in\mathcal{P}(E)$ and an operator $A\in
\textbf{E}nd(E)$.

(ii) For $A, B\in \textbf{E}nd(E)$ and $\mu ,\nu\in\mathcal{P}(E)$
we have equality
\[
\mathcal{B}^{(\nu)} \cdot \mathcal{A}^{(\mu)} =
\mathcal{C}^{(\mu\ast A \nu)}, \ \mbox{where}\ C:=AB \  ( "\cdot "\
\mbox{means the composition of operators}. )
\]

(iii) For one-parameters families of operators $A_t\in
\textbf{E}nd(E)$ and probability measures $\rho_t\in \mathcal{P}(E)$
($t\ge0)$,
\[
[\,\mathcal{A}_s^{(\rho_s)} \cdot  \mathcal{A}_t^{(\rho_t)}
=\mathcal{A}_{t+s}^{(\rho_{t+s})}\,]\ \ \mbox{iff}\ \ [\,A_t \cdot
A_s=A_{t+s}, \ \ \mbox{and}\ \ \rho_{t+s}=\rho_t\ast A_t\rho_s\,]
\]
\end{prop}
\emph{Proof.} (i). Suppose
$\mathcal{A}^{(\mu)}=\mathcal{B}^{(\nu)}$. Putting $x =0$ in (8) we
have
\[
(\mathcal{A}^{(\mu)}f)(0)=\int_E\,f(y)\mu(dy)= \int_E\,f(y)\nu(dy),
\ \ \mbox{for all}\ \ f\in C_b(E)
\]
which, by Riesz Theorem, implies that $\mu=\nu$. Furthermore, since
\[
(\mathcal{A}^{(\mu)}f)(x)=\int_E\int_Ef(u+y)\delta_{Ax}(du)\mu(dy)=\int_E
f(z)(\delta_{Ax}\ast\mu)(dz),
\]
we conclude that $\delta_{Ax}\ast\mu=\delta_{Bx}\ast\mu$, for all
$x\in E$, that is, $A=B$.

(ii). For $f\in C_b(E)$ and $x \in E$, by (8), we have
\begin{multline*}
((\mathcal{B}^{(\nu)} \cdot \mathcal{A}^{(\mu)})f)(x)=
(\mathcal{B}^{(\nu)}(\mathcal{A}^{(\mu)}
f))(x)=\int_E(\mathcal{A}^{(\mu)} f)(Bx+y)\nu(dy)=\\
\int_E\,(\int_E \,f(A(Bx+y)+z))\mu(dz))\nu(dy)= \int_E \int_E
f((AB)x+Ay+z)\mu(dz)\nu(dy)\\= \int_Ef((AB)x+u)(\mu \ast A\nu)(du)=
(\mathcal{C}^{(\mu \ast A\nu)}f)(x),
\end{multline*}
which proves (ii).

(iii) Since, by (ii),  \  $\mathcal{A}_s^{(\rho_s)} \cdot
\mathcal{A}_t^{(\rho_t)}= (\mathcal{A}_t\mathcal{A}_s)^{(\rho_t\ast
A_t\rho_s)}$ \ \ thus in order to have the equality in (iii) it is
necessary and sufficient that $A_tA_s=A_{t+s}$ and $\rho_t\ast A_t
\rho_s=\rho_{t+s}$, because of (i).

\medskip
For  a given $C_0$-semigroup $(T_t, t\ge0)$ on $E$ and a family of
probability measures $\rho_t$,  one-parameter semigroups
$\mathcal{T}_t\equiv \mathcal{T}_t^{(\rho_t)}$ (on $C_b(E)$) are
called \emph{one-parameter generalized Mehler semigroup}. Hence
necessarily and sufficiently one has: \ $\rho_{t+s}=\rho_t\ast
T_t\rho_s=\rho_s\ast T_s\rho_t$ \ for all \ $t,s \ge0$; comp.
Proposition 3(iii). Such equations were called \emph{cocycles} in
Hofmann-Jurek (1996)).

Explicitly we can write,
\begin{equation}
(\mathcal{T}_t\,f)(x)=\int_E\,f(T_tx+y)\rho_t(dy), \qquad f \in
C_b(E) .
\end{equation}
Firstly, let us note that we have the following relation:
\begin{cor}
Each Urbanik semigroup $\textbf{D}(\mu)$, $\hat{\mu}(y)\neq 0, y \in
E^{\prime}$, that contains one-parameter $C_0$-semigroup $(T_t,
t\ge0)$ induces a generalized Mehler semigroup $(\mathcal{T}_t,
t\ge0)$ by taking in (9) $\rho_t=\nu_{T_t}$ from (4).
\end{cor}
Secondly, inspired by the technique from the operator-limit
distribution theory we get

\begin{prop} If \ $t \to \rho_t$ is continuous at  zero and
$\rho_{t+s}=\rho_t\ast T_t\rho_s, \ \mbox{for all} \ t,s \ge 0$ \
(cocycle equation) then there exist a cadlag process \ $Z(t),
t\ge0$, \ with independent increments such that \
$\mathcal{L}(Z(t))=\rho_t$ \ and \ $Z(0)=0$ a.s. In particular, all
$\rho_t$ are infinitely divisible.
\end{prop}

\emph{Proof.} From the Kolmogorov's Extension Theorem, (on a family
of consistent distributions), in order to describe (in distribution)
a process $Z_t, t\ge0,$ starting from zero (i.e. $Z_0=0$ with P.1)
and with independent increments it is necessary and sufficient to
give the probability distributions of all increment $(Z_t-Z_s,\ t\ge
s\ge0)$ (in particular, one gets distributions of $Z_t$) in a such
way that $Z_t\stackrel{d}{=} (Z_t-Z_s) +Z_s$ with the the summands
independent for all $t\ge s\ge 0$.

Let us define
\begin{equation}
\mathcal{L}(Z_t-Z_s):=T_s\rho_{t-s}, \ \ t\ge s\ge 0, \ \mbox{in
particular} \ \ Z_t\stackrel{d}{=}\rho_t.
\end{equation}
Then, by the independence and the cocycle equation we get

$(Z_s-Z_0)+(Z_t-Z_s)\stackrel{d}{=}\rho_s\ast T_s\rho_{t-s}
=\rho_t\stackrel{d}{=}Z_t$.

\noindent Since $t\to\rho_t$ is continuous therefore Banach
space-valued process $Z_t$ is continuous in probability.
Consequently, for $Z_t, t\ge0$, there exist its c\'adl\'ag version
(in French \emph{c\'adl\'ag$\equiv$ \underline{c}ontinu
\underline{\'a} \underline{d}roite avec des \underline{l}imites
\underline{\'a} \underline{g}auche}, i.e., paths are right
continuous with left-hand limits) $Z(t),t\ge0$; cf. Jurek-Vervaat
(1983), Theorem A.1.1 on p. 260.

\noindent Finally, using (10) for each $t\ge0$ and each $n\ge1$ we
have
\[
\rho_t\stackrel{d}{=}Z(t)=\sum_{k=1}^n \big(Z(\frac{k\,t}{n})-
Z(\frac{(k-1)t}{n})\big)\stackrel{d}{=}\rho_{\frac{t}{n}}\ast
T_{t/n}\,\rho_{\frac{t}{n}}\ast T_{2t/n}\,\rho_{\frac{t}{n}}...\ast
T_{(n-1)t/n}\,\rho_{\frac{t}{n}},
\]
and the triangular array is infinitesimal, i.e, for each
$\epsilon>0$
\[
\lim_{n\to\infty}\,\max_{0\le j\le
n-1}\,\,(T_{jt/n}\,\rho_{\frac{t}{n}})(||x||\ge \epsilon)=0,
\]
because of  (2) and the fact that $\rho_s\Rightarrow \delta_0$, as
$s\to 0$. This proves the infinite divisibility of $\rho_t$ and thus
completes the proof.

\begin{rem}
 Note that in Proposition 4 the infinite divisibility one gets
 from the stochastic independence of increments of the process $Z(t), t\ge0$ as well.
 The infinite divisibility  in the cocycle equations was proved in
 Schmuland-Sun (2001) by different (analytic) methods and without the continuity condition. In our
 approach the continuity $t\to \rho_t$ was used to get c\'adl\'ag paths
 of the constructed process and consequently the infinite divisibility property.
\end{rem}

\begin{cor}
Each generalized Mehler semigroup $\mathcal{T}_t$, with \
$t\to\rho_t$ continuous at zero, is of the form
\[
(\mathcal{T}_t\,f)(x)=\textbf{E}\,[\,f(T_tx+Z(t))\,] \ \ f\in
C_b(E),
\]
for some $C_0$-semigroup $(T_t, t\ge0)$ in \textbf{E}nd(E) and some
E-valued c\'adl\'ag process $(Z(t), t\ge0)$ with independent
increments, $Z(0)=0$ a.s. and  \
$Z(t)-Z(s)\stackrel{d}{=}T_s\,Z(t-s)$ for all $t\ge s\ge 0$.
\end{cor}

\medskip
For an \textbf{E}nd(E)- valued function $g(t)$ of locally bounded
variation and E-valued c\'adl\'ag process with independent
increments $(Y(t),t\ge0)$, let us define \emph{a random integral}
via formal integration by parts formula:
\begin{equation}
\int_{(a,b]}g(t)d\,Y(t):= g(b)Y(b)-g(a)Y(a)-\int_{(a,b]}dg(t)Y(t),
\end{equation}
where the right-hand side is defined as path-wise approximation by
partial sums of the form $\sum_{j=1}^n\,(g(t_j)-g(t_{j-1}))Y(t_j)$;
in a similar way as in Jurek (1982) and Jurek-Vervaat (1983) or
 Jurek-Mason (1993).

\medskip
Recall that by stochastic \emph{L\'evy process} we mean c\'adl\'ag
process with stationary and independent increment, and starting from
zero.
\begin{thm}
For each $C_0$-semigroup $\textbf{T}=(T_t,t\ge0)$ and each E-valued
L\'evy process $\textbf{Y}=(Y(t),t\ge0)$ a semigroup of operators
\begin{equation}
\mathcal{T}^{\textbf{T,Y}}f\equiv(\mathcal{T}_t\,f)(x):=\mathbf{E}
[f(T_tx+\int_{(0,t]}T_{t-s}\,dY(s))], \ f\in C_b(E),
\end{equation}
is a generalized Mehler semigroup.

Conversely, if $\mathcal{T}$ is a generalized Mehler semigroup with
the mapping $t\to \rho_t$ continuous  at zero and an one-parameter
\underline{\emph{group}} $\textbf{T}=(T_t, t \in \Rset)$ of
operators then there exist a unique (in distribution) L\'evy
c\'adl\'ag process $Y(t), t\ge0$ such that
$\mathcal{T}=\mathcal{T}^{\textbf{T,Y}}$.
\end{thm}
\emph{Proof.} Let $V_t$ denotes the random integral part in (12) and
let  $\rho_t$ be its probability distribution. Then, using the
standard argument of approximation by partial sums, we infer that
\begin{multline}
\log \hat{\rho_t}(y):=\log\,\mathbf{E}[e^{
i<y,V_t>}]=\int_{(0,t]}\log\,\mathbf{E}[e^{i<T^*_{t-s}y,Y(1)>}]ds \\
= \int_{(0,t]}\log\mathbf{E}[e^{i<T_{r}^*\,y,Y(1)>}]dr. \qquad
\qquad
\end{multline}
Hence by a simple calculations we get
$\log\hat{\rho_t}(y)+\log\hat{\rho_s}(T^*_t\,y)=
\log\hat{\rho}_{t+s}(y)$, $y\in E^{\prime}$. Thus the family
$\rho_t, t\ge 0$ satisfies the cocycle equation and therefore (12)
defines a Mehler semigroup.

Conversely, let a generalized Mehler semigroup $\mathcal{T}$ be
given by (9). Then by Corollary 2, there exists a c\'adl\'ag process
$(Z(t),t\ge0)$ with independent increments. Furthermore, the
stochastic process
\begin{equation}
Y(t):=\int_{(0,t]}T_{-s}\,dZ(s), t\ge 0
\end{equation}
is with independent increments , because so is $Z(.)$. And more
importantly, for $t>s$, using the fact that
$T_{-s}(Z(t)-Z(s))\stackrel{d}{=}Z(t-s)$, by (10), we conclude
\[
Y(t)-Y(s)=\int_{(0,t-s]}T_{-v}T_{-s}dZ(v+s)\stackrel{d}{=}\int_{(0,t-s]}T_{-v}dZ(v)=Y(t-s),
\]
i.e. that $Y$ is a L\'evy process (independent and stationary
increments). Since, by (13) and (14), we have that
\[
\int_{(0,\,t]}
T_{t-s}dY(s)\stackrel{d}{=}\int_{(0,\,t]}T_sdY(s)=Z(t), \ \ \
\mbox{for each $t>0$},
\]
which with Corollary 2 gives the formula (12). The uniqueness in
distribution is a consequence of the Kolmogorov's Extension Theorem
and therefore the proof is complete.

\begin{rem}
Note that for a $C_0$-semigroup $T_t$, on Banach space $E$, and for
an E-valued L\'evy c\'adl\'ag stochastic process $Y$, the processes
given by random integrals
\[
V(t):=\int_{(0,\,t]} T_{t-s}dY(s), \ \ \
Z(t):=\int_{(0,\,t]}T_sdY(s), \ \ t\ge 0,
\]
have only identical marginal (one-dimensional) distributions. The
process $Z$ has independent increments while $V$ is a Markov
process.
\end{rem}

\begin{cor}
(a)\ On Euclidean spaces ($E=\Rset^d$) all generalized Mehler
semigroups are of the form (12).

(b)\ On arbitrary separable Banach space, for any uniformly
continuous semigroups $T_t=\exp tQ$ ($Q$ bounded operator), all
generalized Mehler semigroups are of the form (12).

\end{cor}
\begin{rem}
(a) In the above proof of Theorem 1, the group property of $(T_t, t
\in \Rset)$ was only used to define  the L\'evy process \ $Y$ \ (in
(14)) via the additive process $Z$ from Corollary 2. However, what
we only need is that, for the given additive process $Z$ (in
Corollary 2), the stochastic differential equation
\[
\ \  T_t\,dX(t)=dZ(t),\ \ X(0)=0,    \ \ \  \  \ \ (*)
\]
has a solution for $X(\cdot)$ in L\'evy processes. But since not all
additive processes are semimartingles, the stochastic equations of
the form (*) may be not solvable; cf. Jacod and Shiryaev (1987),
Theorems 4.14 and 4.15 p. 106.

(b) The representation like the above (12) can be derived from
Bogachev, Roeckner and Schmuland (1996), Lemma 2.6  but only for
Hilbert spaces $H$ and more importantly, under restrictive
assumptions that, for $y\in H$, the functions $
t\to\widehat{\rho_t}(y)$ are differentiable at zero. The same
setting is also in Fuhrman and Roeckner (2004).
\end{rem}

Let us also recall here that stochastic processes, from Theorem 1,
\[
\ \ \ \ \  U_t:= T_tx+\int_{(0,t]}T_{t-s}\,dY(s), \ \  t\ge0,
\]
are called \emph{Ornstein-Uhlenbeck processes} and are well studied
on Hilbert spaces. [These are solutions to so called Langevin
equations.] Furthermore, one can easily express L\'evy-Khintchine
formula of $U_t$ in terms of the corresponding parameters of $Y(1)$;
cf. for instance Jurek-Mason (1993), Section 3.6 in Chapter 3.

\medskip
\textbf{Open problem.} It is known that on arbitrary separable
Banach space if $\lim_{t\to \infty}\exp\,tV =0$ (in the operator
topology; $V$ is a bounded operator) then, for a c\'adl\'ag L\'evy
process $Y$, the limit $\lim_{t\to \infty}\int_{(0,t]}e^{sV}dY(s)$
exists (almost surely, in probability or in distribution) if and
only if $\textbf{E}[\log(1+||Y(1)||]<\infty$, cf. Jurek (1982).

Is there an analogous criterium true for $C_0$-semigroups \textbf{T}
on $E$ ? Or how does the existence of a limit depend on the
infinitesimal generator \emph{J} of  the semigroup \textbf{T} and
the variable $Y(1)$?

Recall that Applebaum (2005) in Theorem 9, showed that
log-integrability is sufficient for \emph{exponentially stable
contraction semigroup $(T_t, t\ge 0)$) on a Hilbert space}. His
proof uses the L\'evy-It\^o  decomposition of the process $Y$.

\medskip
\textbf{Acknowledgements.} Professor Jerzy Zabczyk (Polish Academy
of Sciences, Warsaw) called my attention to generalized Mehler
semigroups. In December 2004  I gave a lecture at his seminar and
for that occasion I prepared lecture notes \emph{Measure valued
cocycles from my papers in 1982 and 1983 and Mehler semigropus}; cf.
www.math. uni.wroc.pl/$^{\sim}$zjjurek . At Zabczyk's suggestion I
mailed those notes to Professor M. Roeckner (Bielefeld) who in turn
mailed me some of his papers on Mehler semigropus. Later on, I also
had some useful discussions with Professor D. Applebaum (Sheffield).

\medskip
\medskip
\begin{center}
\textbf{REFERENCES}
\end{center}
\noindent[1] D. Applebaum (2005). Martingale-valued measures,
Ornstein-Uhlenbeck processes with jumps and operator
self-decomposability in Hilbert spaces, in: \emph{S\'eminaire de
Probabilit\'es} vol. 39 , Lecture Notes in Math. no 1874, Springer,
pp. 171-196.

\noindent [2] A. Araujo and E. Gine (1980). \emph{The central limit
theorem for real and Banach valued random variables.} John Wiley \&
Sons, New York.

\noindent [3] V. I. Bogachev, M. Roeckner and B. Schmuland (1996).
Generalized Mehler semigroups and applications, \emph{Probab. Th.
Rel. Fields} vol. 105, pp. 193-225.

\noindent[4] A. Chojnowska-Michalik (1985). Stationary distributions
for $\infty$-dimenional linear equations with general noise, in:
\emph{Lecture Notes on Control and Inf. Sci.}, vol. 69,
Springer-Verlag, Berlin , pp.14-25.

\noindent[5] A. Chojnowska-Michalik (1987). On processes of
Ornstein-Uhlenbeck type in Hilbert space, \emph{Stochastics} vol.
21, pp. 252-286.

\noindent [6] M. Fuhrman and M. Roeckner (2000). Generalized Mehler
semigroups: the non-Gaussian case, \emph{Potential Analysis}
\textbf{12}, pp.1-47.

\noindent [7] K. H. Hofmann and Z. J. Jurek (1996). Some analytic
semigroups occurring in probability theory, \emph{J. Theoretical
Probab. 9}, No. 3, pp. 745-763.

\noindent[8] J. Jacod and A. N. Shiryaev (1987).\emph{Limit theorems
for stochastic processes}, Springer-Verlag ,  Berlin.

\noindent[9] Z. J. Jurek (1982).  An integral representation of
operator-selfdecomposable random variables, \emph{ Bull. Acad. Pol.
Sci; Serie Math.} XXX no.7-8, pp. 385-393.

\noindent[10] Z. J. Jurek (1983). Limit distributions and
one-parameter groups of linear operators on Banach spaces, \emph{J.
Multivariate Analysis}, 13, no. 4, pp. 578-604.

\noindent[11] Z. J. Jurek (1988). Random integral representations
for classes of limit distributions similar to L\'evy class $L_0$,
\emph{Probab. Th. Fields}, 78, pp. 473-490.

\noindent[12] Z. J. Jurek (2004). Measure valued cocycles from my
papers in 1982 and 1983 and Mehler semigroups;
www.math.uni.wroc.pl/$^{\sim}$zjjurek

\noindent[13] Z. J. Jurek and J. D. Mason (1993).
\emph{Operator-limit distributions in probability theory.} John
Wiley \&Sons, New York.

\noindent [14] Z. J. Jurek and W. Vervaat (1983). An integral
representation for selfdecomposable Banach space valued random
variables, \emph{Z. Wahrscheinlichkeitstheorie verw. Gebiete}, 62,
pp. 247-262.

\noindent[15] Z. H. Li (2002). Skew convolution semigroups and
related immigration processes, \emph{Theory Probab. Appl.} 46, pp.
274-296.

\noindent [16] W. Linde (1986). \emph{Probability in Banach
spaces--stable and infinite divisible distributions}, Wiley, New
York.

\noindent[17] K. R. Parthasarathy (1967).\emph{Probability measures
on metric spaces}. Academic Press, New York and London.

\noindent [18] B. Schmuland and W. Sun (2001). On equation
$\mu_{t+s}=\mu_s\ast T_s\mu_t$,\ \emph{Stat. Probab. Letters},
vol.52, 183-188.

\noindent[19] K. Urbanik (1972). L\'evy's probability measures on
Euclidean spaces, \emph{Studia Math.} 44, pp. 119-148.

\noindent[20] K. Urbanik (1978). L\'evy's probability measures on
Banach spaces, \emph{Studia Math.} 63, pp. 283-308.

\medskip
\noindent
Institute of Mathematics \\
University of Wroc\l aw \\
Pl.Grunwaldzki 2/4 \\
50-384 Wroc\l aw, Poland \\
e-mail: zjjurek@math.uni.wroc.pl\\
www.math.uni.wroc.pl/$^{\sim}$zjjurek

\end{document}